\documentclass{amsart}
\usepackage{amsfonts}
\usepackage{amsmath}
\usepackage{enumerate}
\usepackage{amsthm}
\usepackage{hyperref}
\usepackage{cite}

\begin{document}

\newtheorem{definition}{Definition}[section]
\newtheorem{theorem}[definition]{Theorem}
\newtheorem{proposition}[definition]{Proposition}
\newtheorem{remark}[definition]{Remark}
\newtheorem{lemma}[definition]{Lemma}
\newtheorem{corollary}[definition]{Corollary}
\newtheorem{example}[definition]{Example}

\numberwithin{equation}{section}

\title[Prescribed Chern scalar curvatures...]{Prescribed Chern scalar curvatures on complete Hermitian manifolds}
\author[W. Yu]{Weike Yu}
\date{}
\thanks{The author is supported by NSFC (No. 12501075), Basic Research Program of Jiangsu (No. BK20250644), Natural Science Foundation of the Jiangsu Higher Education Institutions of China (No. 25KJB110008), the start-up research funds from Nanjing Normal University with account No. 184080H201B160, and the Laboratory of Mathematics for Nonlinear Science, Fudan University.}

\begin{abstract}
In this paper, we investigate the problem of prescribing the Chern scalar curvature on complete noncompact Hermitian manifolds, and generalize the existence results of Aviles--McOwen [J. Differential Geom., 21 (1985): 269-281] from Poincar\'e disks to higher dimensional Hermitian manifolds.
\end{abstract}
\maketitle
\section{Introduction}
Let $(M^n, J, h)$ be a complex manifold of complex dimension $n\geq2$ with a Hermitian metric $h$ and a complex structure $J$. The fundamental form $\omega$ is defined by $\omega(\cdot, \cdot)=h(J\cdot, \cdot)$. In this paper, we will not distinguish between the Hermitian metric $h$ and its corresponding fundamental form $\omega$. On a Hermitian manifold $(M, \omega)$, there exists a unique linear connection $\nabla^{Ch}$ which is called the Chern connection preserving both the Hermitian metric $h$ and the complex structure $J$, whose torsion $T^{Ch}$ has vanishing $(1,1)$-part everywhere. The scalar curvature with respect to the Chern connection (referred to as the Chern scalar curvature) can be given by
\begin{align}
S^{Ch}(\omega)=\text{tr}_{\omega}Ric^{(1)}(\omega)=\text{tr}_\omega \sqrt{-1}\partial \bar{\partial} \log \omega^n,
\end{align}
where $Ric^{(1)}(\omega)$ denotes the first Chern Ricci curvature, and $\omega^n$ is the volume form of $(M^n, \omega)$.

A basic problem in differential geometry is that of studying the set of curvature functions that a manifold possesses. The following problem is a Hermitian analogue of prescribing scalar curvatures:
\begin{center}
\begin{minipage}{8cm}
If we are given a smooth real-valued function $f$ on $(M, \omega)$, is there a Hermitian metric $\tilde{\omega}$ conformal to $\omega$ such that its Chern scalar curvature $S^{Ch}(\tilde{\omega})=f$?
\end{minipage}
\end{center}
When $f$ is a constant, the above problem is referred to as the Chern-Yamabe problem,  which was first proposed by Angella--Calamai--Spotti in \cite{[ACS]}, and they proved that it is solvable on any compact Hermitian manifolds with nonpositive Gauduchon degree (see \eqref{2.6...}). If $f$ is a general function, the above problem for compact Hermitian manifolds has been studied by several authors (see, e.g., \cite{[Bar], [CZ], [Fus], [Ho], [LM], [Yu1], [Yu2], [Yu3], [Yu4]}). Besides, the problem has also been extended to the almost Hermitian case (cf. \cite{[LU], [LZZ]}) and to a special class of noncompact manifolds under certain analytical conditions (cf. \cite{[WZ1], [WZ2]}). Note that if $M$ is a Riemannian surface (i.e., $\dim_{\mathbb{C}}M=1$), the above problem reduces to the prescribed Gaussian curvature problem, which has been extensively studied (see, e.g., \cite{[AM-1], [Aub], [Ber], [BEW], [BGS], [CP], [DL], [KW1], [KW2], [KW3], [Ma], [Mor], [PR]}).

In the present paper, we will investigate the problem of prescribing the Chern scalar curvature on complete noncompact Hermitian manifolds $(M^n, \omega)$, which is equivalent to solving the following partial differential equation:
\begin{align}\label{1.2}
-\Delta^{Ch}_\omega u+S^{Ch}(\omega)=S^{Ch}(\tilde{\omega})e^{\frac{2}{n}u},
\end{align}
where $\tilde{\omega}=e^{\frac{2}{n}u}\omega$ for some $u\in C^\infty(M)$, and $\Delta^{Ch}_\omega$ denotes the Chern Laplacian with respect to $\omega$. Roughly speaking, we will establish some existence, uniqueness and nonexistence results for this problem. 

On a complete noncompact Hermitian manifold $(M^n, \omega)$, we consider the existence of smooth solutions to \eqref{1.2} that are bounded on $M$. This guarantees that the conformal Hermitian metrics $\tilde{\omega}$ and $\omega$ are uniformly equivalent, i.e., $C^{-1}\omega\leq \tilde{\omega}\leq C\omega$ for some constant $C>0$. In particular, $\tilde{\omega}$ is also complete. Using the method of upper and lower solutions and the Omori--Yau maximum principle on complete Hermitian manifolds, we establish the following theorem, which extends the results of Aviles--McOwen \cite{[AM-1]} from Poincar\'e disks to higher dimensional Hermitian manifolds.

\begin{theorem}\label{theorem1.3}
Let $(M^n, \omega)$ be a complete noncompact Hermitian manifold with its Chern scalar curvature $S^{Ch}(\omega)$ satisfying
\begin{align}
-C_1\leq S^{Ch}(\omega)\leq -C_2,
\end{align}
where $C_1, C_2$ are two positive constants. Suppose that $K$ is a smooth function on $M$ with 
\begin{equation}
\begin{aligned}
&K(x)\leq 0, \quad \forall\ x\in M,\\
-a^2\leq  &K(x)\leq -b^2, \quad\forall\ x\in M\setminus D,
\end{aligned}
\end{equation}
where $D$ is a compact subset of $M$, $a\geq b>0$ are two constants. Then there exists a complete Hermitian metric $\tilde{\omega}$ which is conformal and uniformly equivalent to $\omega$, such that its Chern scalar curvature $S^{Ch}(\tilde{\omega})=K$. If, in addition, we assume that $(M^n, \omega)$ satisfies the second Chern Ricci curvature
\begin{align}\label{1.5.}
Ric^{(2)}(X, \overline{X})\geq-C_3(1+r(x))^2,
\end{align}
and the torsion 
\begin{align}\label{1.6.}
\|T^{Ch}(X,Y)\|\leq C_4(1+r(x)),
\end{align}
where $X, Y\in T^{1,0}_{x}M$ satisfy $\|X\|=\|Y\|=1$, $C_3, C_4$ are two positive constants, $r(x)$ is the Riemannian distance between a fixed point $x_0$ and $x$ in $(M, \omega)$, then the above metric $\tilde{\omega}$ is unique among all conformal metrics that are uniformly equivalent to $\omega$.
\end{theorem}

When $\dim_{\mathbb{C}}M=1$, the Chern scalar curvature is equal to twice the Gaussian curvature. Hence, we recover the Aviles--McOwen's result on Poincar\'e disks.

\begin{corollary}[\cite{[AM-1]}]
Let $M=\{x\in \mathbb{R}^2: |x|<1\}$ with the metric $g_{ij}=4(1-|x|^2)^{-2}\delta_{ij}$ which has constant Gaussian curvature $k\equiv-1$. Suppose $K$ is a smooth function on $M$ with 
\begin{equation}
\begin{aligned}
&K(x)\leq 0, \quad \forall\ x\in M,\\
-a^2\leq & K(x)\leq -b^2, \quad\forall\ x\in M\setminus D,
\end{aligned}
\end{equation}
where $D$ is a compact subset of $M$, $a\geq b>0$ are two constants. Then there exists a unique complete K\"ahler metric $\tilde{g}$ which is conformal and uniformly equivalent to $g$, such that its Gaussian curvature is $K$.
\end{corollary}

\begin{remark}
In \cite{[AM-1]}, Aviles and McOwen proved the above corollary by using the method of upper and lower solutions. Indeed, they constructed a lower solution $u_-=-\log a$ and an upper solution $u_+=\alpha(1-|x|^2)-\log b$, where one may assume $b\in (0,1)$ and choose $\alpha>0$ such that $\alpha\geq \max_{x\in D}\frac{K(x)+1}{(1-|x|^2)^{2}}$. Unfortunately, their construction of the upper solution $u_+$ is not directly applicable to the higher dimensional case, as it relies on the Euclidean distance in the ambient space $\mathbb{R}^2$. In the present paper, we construct an upper solution and a lower solution on general higher dimensional Hermitian manifolds, and then apply the method of upper and lower solutions to establish the existence result.
\end{remark}

Making use of Theorem \ref{theorem1.3}, we also establish the following existence and uniqueness result for metrics of constant negative Chern scalar curvature.
\begin{theorem}\label{theorem1.6}
Let $(M^n, \omega)$ be a complete noncompact Hermitian manifold with its Chern scalar curvature $S^{Ch}(\omega)$ satisfying
\begin{equation}\label{1.6}
\begin{aligned}
&S^{Ch}(\omega)(x)\leq 0, \quad\forall\ x\in M\\
-a^2\leq &S^{Ch}(\omega)(x)\leq -b^2, \quad\forall\ x\in M\setminus D,
\end{aligned}
\end{equation}
where $D$ is a compact subset of $M$, $a\geq b>0$ are two constants. Then there exists a complete Hermitian metric $\tilde{\omega}$ which is conformal and uniformly equivalent to $\omega$, such that its Chern scalar curvature $S^{Ch}(\tilde{\omega})\equiv-1$. If, in addition, we assume that $(M, \omega)$ satisfies \eqref{1.5.} and \eqref{1.6.}, then $\tilde{\omega}$ is unique among all conformal metrics that are uniformly equivalent to $\omega$.
\end{theorem}

If we remove the assumption $-a^2\leq S^{Ch}(\omega)$ in \eqref{1.6}, we can still establish an existence result by following Aviles--McOwen's idea in \cite{[AM-2]}, although the equations studied are different. However, the Hermitian metric $\tilde{\omega}$ found in the conformal class of $\omega$ may no longer be uniformly equivalent to $\omega$; it still satisfies $\tilde{\omega}=e^{\frac{2}{n}u}\omega\geq C\omega$ for some constant $C>0$. This leads to the conclusion that $\tilde{\omega}$ is also complete.

\begin{theorem}\label{theorem1.7.}
Let $(M^n, \omega)$ be a complete noncompact Hermitian manifold with its Chern scalar curvature $S^{Ch}(\omega)$ satisfying
\begin{equation}
\begin{aligned}
&S^{Ch}(\omega)(x)\leq 0, \quad\forall\ x\in M\\
&S^{Ch}(\omega)(x)\leq -b^2, \quad\forall\ x\in M\setminus D,
\end{aligned}
\end{equation}
where $D$ is a compact subset of $M$, $b>0$ is a constant. Then there exists a complete Hermitian metric $\tilde{\omega}$ which is conformal to $\omega$ and satisfies $\tilde\omega\geq C\omega$ for some constant $C>0$, such that its Chern scalar curvature $S^{Ch}(\tilde{\omega})\equiv-1$.
\end{theorem}
\begin{remark}
In \cite{[Yuan]}, the author studied fully nonlinear elliptic equations on Hermitian manifolds and applied these results to obtain geometric conclusions about the existence of conformal Hermitian metrics with prescribed Chern-Ricci curvature functions. As a special case, the above theorem was also obtained. However, their method is completely different from that used in this paper.
\end{remark}

Finally, by applying the generalized maximum principle, we obtain a nonexistence result for prescribed Chern scalar curvature on complete noncompact Hermitian manifolds.

\begin{theorem}\label{theorem1.9}
Let $(M^n, \omega)$ be a complete noncompact Hermitian manifold with 
\begin{align}
Ric^{(2)}(X, \overline{X})\geq-C_1(1+r(x))^2,
\end{align}
\begin{align}
S^{Ch}(\omega)\geq 0,
\end{align}
and the torsion 
\begin{align}
\|T^{Ch}(X,Y)\|\leq C_2(1+r(x)),
\end{align}
where $X, Y\in T^{1,0}_{x}M$ satisfy $\|X\|=\|Y\|=1$, $r(x)$ is the Riemannian distance between a fixed point $x_0$ and $x$ in $(M, \omega)$, and $C_1, C_2$ are two positive constants. Let $K$ be a smooth function on $M$ satisfying
\begin{equation}
\begin{aligned}
&K(x)\leq 0, \quad \forall\ x\in M,\\
&K(x)\leq -b^2, \quad\forall\ x\in M\setminus D,
\end{aligned}
\end{equation}
where $D$ is a compact subset of $M$, $b>0$ is a constant. Then the Hermitian metric $\omega$ cannot be conformally deformed to another Hermitian metric $\tilde{\omega}$ such that its Chern scalar curvature is $K$.
\end{theorem}

This paper is organized as follows. In Section 2, we recall some basic concepts and notation related to the prescribed Chern scalar curvature problem. In Section 3, we study the prescribed Chern scalar curvatures on complete noncompact Hermitian manifolds, and prove Theorem \ref{theorem1.3}, Theorem \ref{theorem1.6} and Theorem \ref{theorem1.7.}. In Section 4, we establish a nonexistence result (Theorem \ref{theorem1.9}) for this problem on complete noncompact Hermitian manifolds.

\textbf{Acknowledgments.} The author would like to thank Prof. Yuxin Dong and Prof. Xi Zhang for their continued support and encouragement.

\section{Preliminaries}

Suppose that $(M^n, J, h)$ is a Hermitian manifold with complex dimension $n$ and its fundamental form $\omega(\cdot,\cdot)=h(J\cdot, \cdot)$. Let $TM^{\mathbb{C}}=TM\otimes \mathbb{C}$ be the complexified tangent space of $M$, and we extend $J$ and $h$ from $TM$ to $TM^{\mathbb{C}}$ by $\mathbb{C}$-linearity. Then we have the decomposition
\begin{align}
TM^{\mathbb{C}}=T^{1,0}M\oplus T^{0,1}M,
\end{align}
where $T^{1,0}M$ and $T^{0,1}M$ are the eigenspaces of the complex structure $J$ corresponding to the eigenvalues $\sqrt{-1}$ and $-\sqrt{-1}$, respectively. Furthermore, every $m$-form can also be decomposed into $(p, q)$-forms for all $p, q \geq 0$ with $p + q = m$ by extending the complex structure $J$ to act on forms.

On a Hermitian manifold $(M^n, J, h)$, there exists a unique affine connection $\nabla^{Ch}$, called the Chern connection, which preserves both the Hermitian metric $h$ and the complex structure $J$, that is,
\begin{align}
 \nabla^{Ch}h=0,\quad \nabla^{Ch}J=0, 
\end{align}
whose torsion $T^{Ch}(X,Y)=\nabla_XY-\nabla_YX-[X,Y]$ satisfies 
\begin{align}
T^{Ch}(JX,Y)=T^{Ch}(X,JY)
\end{align}
for any $X, Y\in TM$. In this paper, we do not distinguish between the Hermitian metric $h$ and its corresponding fundamental form $\omega$. 

For a Hermitian manifold $(M^n, \omega)$ with Chern connection $\nabla^{Ch}$, following the notation in \cite[Sec. 2]{[Yu0]}, the first Chern Ricci curvature is defined by
\begin{align}
Ric^{(1)}=\sqrt{-1}Ric^{(1)}_{i\bar{j}}\theta^i\wedge\theta^{\bar{j}}
\end{align}
with components 
\begin{align}
Ric^{(1)}_{i\bar{j}}=\sum_{k}R_{k\bar{k}i\bar{j}},
\end{align}
where $R_{i\bar{j}k\bar{l}}=h(R(e_k,e_{\bar{l}})e_i, e_{\bar{j}})$ are the components of curvature tensor $R$ of $\nabla^{Ch}$ under a local unitary frame $\{e_i\}_{i=1}^n$ of $T^{1,0}M$, and $\{\theta^i\}_{i=1}^n$ is its dual frame. It is well-known that $Ric^{(1)}$ can be written as
\begin{align}
Ric^{(1)}=\sqrt{-1}\bar{\partial}\partial \log{\omega^n},
\end{align}
where $\omega^n$ is the volume form of $(M^n, \omega)$. In addition, we define the second Chern Ricci curvature by
\begin{align}
Ric^{(2)}=\sqrt{-1}Ric^{(2)}_{i\bar{j}}\theta^i\wedge\theta^{\bar{j}}
\end{align}
with components 
\begin{align}
Ric^{(2)}_{i\bar{j}}=\sum_{k}R_{i\bar{j}k\bar{k}}.
\end{align}
The Chern scalar curvature of $\omega$ is given by
\begin{align}
S^{Ch}(\omega)=\text{tr}_{\omega}Ric^{(1)}(\omega)=\text{tr}_{\omega}Ric^{(2)}(\omega)=\text{tr}_{\omega}\sqrt{-1}\bar{\partial}\partial \log{\omega^n}.
\end{align}

In Hermitian geometry, there is a canonical elliptic operator called the Chern Laplace operator. For any smooth function $u: M\rightarrow \mathbb{R}$, the Chern Laplacian $\Delta^{Ch}_\omega u$ is defined by
\begin{align}
\Delta^{Ch}_\omega u=-2\sqrt{-1} tr_{\omega}\overline{\partial}\partial u.
\end{align}
\begin{lemma}[cf. \cite{[Gau]}]\label{lemma2.1}
On a Hermitian manifold $(M^n,\omega)$, we have
\begin{align}
-\Delta^{Ch}_\omega u=-\Delta_d u+(du,\theta)_{\omega},
\end{align}
where $\Delta_d u=-d^*du$ is the Hodge--de Rham Laplacian, $\theta$ is the Lee form or torsion $1$-form given by $d\omega^{n-1}=\theta\wedge \omega^{n-1}$, and $(\cdot, \cdot)_\omega$ denotes the inner product on $1$-forms induced by $\omega$. Furthermore, if $\omega$ is balanced, i.e., $\theta\equiv0$, then $\Delta^{Ch}_\omega u=\Delta_d u$ for all $u\in C^\infty(M)$.
\end{lemma}
Let 
\begin{align}
\{\omega\}=\{e^{\frac{2}{n}u}\omega\ |\ u\in C^\infty(M)\}
\end{align}
denote the conformal class of the Hermitian metric $\omega$. In \cite{[Gau1]}, Gauduchon proved the following interesting theorem:
\begin{theorem}\label{theorem2.2}
If $(M^n, \omega)$ is a compact Hermitian manifold with complex dimension $n\geq 2$, then there exists a unique Gauduchon metric $\eta\in \{\omega\}$ (i.e., $d^*\theta=0$) with volume $1$.
\end{theorem}
In terms of the above theorem, one can define an invariant $\Gamma(\{\omega\})$ of the conformal class $\{\omega\}$ which is called the Gauduchon degree:
\begin{align}\label{2.6...}
\Gamma(\{\omega\})=\frac{1}{(n-1)!}\int_M c^{BC}_1(K^{-1}_M)\wedge\eta^{n-1}=\int_M S^{Ch}(\eta)d\mu_\eta,
\end{align}
where $\eta$ is the unique Gauduchon metric in $\{\omega\}$ with volume $1$, $c^{BC}_1(K^{-1}_M)$ is the first Bott-Chern class of anti-canonical line bundle $K^{-1}_M$, and $d\mu_\eta$ denotes the volume form of the Gauduchon metric $\eta$. Regarding the sign of \eqref{2.6...}, X. Yang \cite{[Yang]} established the following result:
\begin{theorem}
Let $W$ be the space of Gauduchon metrics on a compact complex manifold $M^n$, and define $F: W\rightarrow \mathbb{R}$ by
\begin{align}
F(\omega)=\int_M S^{Ch}(\omega)\omega^n,\quad \omega\in W.
\end{align}
Let $K_M$ be the canonical line bundle of $M$. Then the image $F(W)$ is one of the following:
\begin{enumerate}
\item $F(W)=\mathbb{R}$ if and only if neither $K_M$ nor $K_M^{-1}$ is pseudo-effective.
\item $F(W)=\mathbb{R}^{>0}$ if and only if $K_M^{-1}$ is pseudo-effective but not unitary flat.
\item $F(W)=\mathbb{R}^{<0}$ if and only if $K_M$ is pseudo-effective but not unitary flat.
\item $F(W)=\{0\}$ if and only if $K_M$ is unitary flat.
\end{enumerate}
\end{theorem}

Given a Hermitian manifold $(M^n, \omega)$, we consider the conformal change $\widetilde{\omega}=e^{\frac{2}{n}u}\omega$. By \cite{[Gau]}, the Chern scalar curvatures of $\widetilde{\omega}$ and $\omega$ have the following relationship:
\begin{align}\label{2.7...}
-\Delta^{Ch}_\omega u+S^{Ch}(\omega)=S^{Ch}(\widetilde{\omega})e^{\frac{2}{n}u},
\end{align}
where $S^{Ch}(\omega)$ and $S^{Ch}(\widetilde{\omega})$ denote the Chern scalar curvatures of $\omega$ and $\widetilde\omega$, respectively. 

To resolve the prescribed Chern scalar curvature problem, it is useful to interpret equation \eqref{2.7...} as the Euler-Lagrange equation of some functional. Unfortunately, in general, such a functional does not always exist. In fact, according to \cite[Prop. 5.3]{[ACS]} and \cite[Prop. 2.12]{[Fus]}, we know the following
\begin{proposition}
Equation \eqref{2.7...} can be seen as an Euler-Lagrange equation for standard $L^2$ pairing if and only if $\omega$ is balanced.
\end{proposition}

 \section{Existence and uniqueness results for prescribed Chern scalar curvature}
 In this section, we will consider the prescribed Chern scalar curvature problem on complete noncompact Hermitian manifolds, and establish some existence and uniqueness results.

\begin{proposition}\label{prop4.1.}
Let $(M^n, \omega)$ be a complete Hermitian manifold. Assume that $f(x, r)\in C^\infty(M\times \mathbb{R})$. If there are two functions $u_+, u_-\in C^0(M)\cap W^{1,2}_{loc}(M)$ satisfying
\begin{align}
&\Delta^{Ch}_\omega u_++f(x, u_+)\leq 0\quad \text{in\ $M$},\label{4.1'}\\
&\Delta^{Ch}_\omega u_-+f(x, u_-)\geq 0\quad \text{in\ $M$},\label{4.2'}\\
&u_+\geq u_-\quad \text{in\ $M$},
\end{align}
in the weak sense, then there exists a function $u\in C^\infty(M)$ such that
\begin{equation}\label{4.4''}
\begin{cases}
\Delta^{Ch}_\omega u+f(x, u)= 0\quad \text{in\ $M$}\\
u_-\leq u\leq u_+ \quad \text{in\ $M$}.
\end{cases}
\end{equation}
Here $u_+$ and $u_-$ are called the upper solution and lower solution of \eqref{4.4''}, respectively.
\end{proposition}
\proof Suppose that $M=\cup_{k\geq 1} \Omega_k$, where each $\Omega_k$ is a bounded domain with smooth boundary, and $\Omega_k\subset\subset \Omega_{k+1}$. Firstly, we will construct a solution of \eqref{4.4''} in $\Omega_k$ by using monotone iteration schemes. Set $r_k^-=\min_{\overline{\Omega}_k}u_-$, $r_k^+=\max_{\overline{\Omega}_k}u_+$, $I_k=[r_k^-, r_k^+]$. Since $f(x, r)\in C^\infty(M\times \mathbb{R})$, there exists a sufficiently large constant $\lambda_k\in \mathbb{R}^+$ and a function $f_k(x, r)\in C^\infty(M\times \mathbb{R})$ such that $f(x, r)=f_k(x, r)-\lambda_k r$, and for each $x\in \Omega_k$, the function $f_k(x, r)$ is increasing with respect to $r\in I_k$. We want to look for a function $u_k$ solving
\begin{equation}\label{4.6'}
\begin{cases}
-\Delta^{Ch}_\omega u_k+\lambda_k u_k= \tilde{f}_k(x, u_k)\quad \text{in\ $\Omega_k$,}\\
u_-\leq u_k\leq u_+ \quad \text{in\ $\Omega_k$},
\end{cases}
\end{equation}
where
\begin{equation}
\tilde{f}_k(x, r)=
\begin{cases}
f_k(x, r^-_k)\quad &\text{if}\ r<r^-_k,\\
f_k(x, r)\quad &\text{if}\ r\in I_k,\\
f_k(x, r^+_k)\quad &\text{if}\ r>r^+_k.
\end{cases}
\end{equation}
For this purpose, we consider a sequence $\{v_-^m\}_{m\geq 1}\subset W^{1,2}(\Omega_k)$ defined by
\begin{equation}\label{4.8'}
\begin{cases}
-\Delta^{Ch}_\omega v^m_-+\lambda_k v^m_-=\tilde{f}_k(x, v^{m-1}_-),\quad \text{in}\ \Omega_k,\\
v^m_--\frac{u_-+u_+}{2}\in W^{1,2}_0(\Omega_k),\\
v^0_-=u_-.
\end{cases}
\end{equation}
According to \eqref{4.1'} and \eqref{4.8'}, we deduce that
\begin{equation}\label{4.9'}
\begin{cases}
(\Delta^{Ch}_\omega-\lambda_k\text{id})(v^1_--u_+)\geq f_k(x, u_+)-f_k(x, u_-)\geq 0\quad \text{in}\ \Omega_k,\\
v^1_--u_+\leq 0\quad \text{mod}\ W^{1,2}_0(\Omega_k),
\end{cases}
\end{equation}
where the condition ``$u\leq v\ \text{mod}\ W^{1,2}_0(\Omega)$'' means $u\leq v+w_0$ for some $w_0\in W^{1,2}_0(\Omega)$. Hence, by the weak maximum principle, we have
\begin{align}
v_-^1\leq u_+.
\end{align}
Indeed, set $w=v^1_--u_+$. According to \cite[Lemma 5.12]{[Gri]}, we have $w\leq 0\ \text{mod}\ W^{1,2}_0(\Omega_k)$ if and only if $w^{+}:=\max\{w, 0\}\in W^{1,2}_0(\Omega_k)$. Taking $w^+$ as a test function of \eqref{4.9'} yields
\begin{align}
-\int_{\Omega_k}|\nabla w^+|^2d\mu_\omega-\int_{\Omega_k}(dw^+, \theta)_\omega w^+d\mu_\omega-\lambda_k \int_{\Omega_k}(w^+)^2d\mu_\omega\geq 0,
\end{align}
where $\theta$ is the Lee form or torsion $1$-form of $\omega$ (see Lemma \ref{lemma2.1}). Thus, by choosing $\lambda_k>-\frac{1}{2}\min_{\overline{\Omega}_k}d^*\theta$, we have
\begin{equation}\label{4.12'}
\begin{aligned}
\int_{\Omega_k}|\nabla w^+|^2d\mu_\omega&\leq -\frac{1}{2}\int_{\Omega_k}(w^+)^2 d^*\theta d\mu_\omega-\lambda_k \int_{\Omega_k}(w^+)^2d\mu_\omega\\
&\leq (-\frac{1}{2}\min_{\overline{\Omega}_k}d^*\theta-\lambda_k)\int_{\Omega_k}(w^+)^2d\mu_\omega\\
&\leq 0.
\end{aligned}
\end{equation}
On the other hand, by the Poincar\'e inequality, there is a constant $C=C(\Omega_k)>0$ such that
\begin{align}\label{4.13'}
\int_{\Omega_k}|\nabla w^+|^2d\mu_\omega\geq C \int_{\Omega_k}(w^+)^2d\mu_\omega.
\end{align}
Combining \eqref{4.12'} and \eqref{4.13'}, we obtain $w^+=0$, i.e., $w=v^1_--u_+\leq 0$. Repeating the above procedures gives
\begin{align}\label{4.13'''}
u_-\leq v^1_-\leq v^2_-\leq \cdots\leq u_+\quad \text{in}\ \Omega_k.
\end{align}
Define $\underline{u}_k=\lim_{m\rightarrow +\infty}v^m_-$ in $\Omega_k$. By the interior elliptic $L^p$-estimates for \eqref{4.8'} and using \eqref{4.13'''}, for any $\Omega'\subset\subset \Omega_k$ there exists a constant $C=C(\Omega', \Omega_k, u_+, u_-, f)>0$ such that $\|v^1_-\|_{W^{2, p}(\Omega')}\leq C$ for any $p>1$. Taking $p=2n+1$ and using the Sobolev embedding theorem $W^{2, 2n+1}(\Omega')\subset C^1(\Omega')$ yields $\|v^1_-\|_{C^{1}(\Omega')}\leq C$. According to the Schauder estimates for \eqref{4.8'} and using \eqref{4.13'''}, we have $\|v^2_-\|_{C^{2, \alpha}(\Omega'')}\leq C$ for any $\Omega''\subset\subset\Omega'\subset\subset \Omega_k$ and $\alpha\in (0,1)$. Repeating the above procedures, we obtain $\|v^m_-\|_{C^{2, \alpha}(\Omega)}\leq C$ for any $\Omega\subset\subset \Omega_k$, $\alpha\in (0,1)$ and $m\geq 2$, where $C$ is a positive constant independent of $m$. By the compact embedding theorem $C^{2, \alpha}(\Omega)\subset\subset C^{2}(\Omega)$, we get  $\underline{u}_k$ is a $C^2$-solution of \eqref{4.6'}. Hence, $\underline{u}_k\in C^2(\Omega_k)$ solves
\begin{equation}
\begin{cases}
\Delta^{Ch}_\omega \underline{u}_k+f(x, \underline{u}_k)=0\quad \text{in\ $\Omega_k$,}\\
u_-\leq \underline{u}_k\leq u_+ \quad \text{in\ $\Omega_k$.}
\end{cases}
\end{equation}

Now let us follow the proof of \cite[Theorem 2.10]{[Ni]} and consider the sequence $\{\underline{u}_k\}_{k\geq 4}$ on $\overline{\Omega}_3$. Since $\underline{u}_k\in I_3$ on $\overline{\Omega}_3$ for any $k\geq 4$ and $f(x, r)\in C^\infty(M\times \mathbb{R})$, there exists a constant $C>0$ independent of $k$ such that $\|f(x, \underline{u}_k)\|_{L^\infty(\overline{\Omega}_3)}<C$ and $\|\underline{u}_k\|_{L^\infty(\overline{\Omega}_3)}<C$. Using the interior elliptic $L^p$-estimates, we get $\|\underline{u}_k\|_{W^{2,p}(\Omega_2)}<C$ for any $p>1$. Taking $p=2n+1$ and using the Sobolev embedding theorem $W^{2, 2n+1}(\Omega_2)\subset C^1(\Omega_2)$ yields $\|\underline{u}_k\|_{C^1(\Omega_2)}\leq C$. According to the Schauder estimates, we have $\|\underline{u}_k\|_{C^{2, \alpha}(\Omega_1)}\leq C$ for any $\alpha\in (0,1)$. From the compact embedding theorem $C^{2, \alpha}(\Omega_1)\subset\subset C^{2}(\Omega_1)$, it follows that there is a subsequence $\{\underline{u}_{k1}\}$ of $\{\underline{u}_k\}$ such that $\underline{u}_{k1}$ converges to $u_1$ in $C^2(\Omega_1)$ as $k\rightarrow +\infty$, where $u_1$ is a solution of
\begin{align}
\Delta^{Ch}_\omega u_1+f(x, u_1)=0\quad \text{in}\ \Omega_1.
\end{align}
Repeating the above procedures with $\{\underline{u}_{k1}\}_{k\geq 5}$ on $\overline{\Omega}_4$, we obtain a subsequence $\{\underline{u}_{k2}\}$ of $\{\underline{u}_{k1}\}$ such that  $\underline{u}_{k2}$ converges to $u_2$ in $C^2(\Omega_2)$ as $k\rightarrow +\infty$, and $u_2$ solves
\begin{equation}
\begin{cases}
\Delta^{Ch}_\omega u_2+f(x, u_2)=0\quad \text{in}\ \Omega_2,\\
u_2|_{\Omega_1}=u_1.
\end{cases}
\end{equation}
Inductively, one can obtain a subsequence $\{\underline{u}_{kj}\}$ which converges to $u_j$ in $C^2(\Omega_j)$ as $k\rightarrow +\infty$. Here $u_j$ is a solution of 
\begin{equation}
\begin{cases}
\Delta^{Ch}_\omega u_j+f(x, u_j)=0\quad \text{in}\ \Omega_j,\\
u_j|_{\Omega_{j-1}}=u_{j-1}
\end{cases}
\end{equation}
for any $j\geq 2$. Define 
\begin{align}
u(x)=\lim_{k\rightarrow +\infty}u_{kk}(x)\quad \ x\in M.
\end{align}
Then $u\in C^2(M)$ and it solves
\begin{equation}
\begin{cases}
\Delta^{Ch}_\omega u+f(x, u)=0\quad \text{in\ $M$}\\
u_-\leq u\leq u_+ \quad \text{in\ $M$},
\end{cases}
\end{equation}
Using the Schauder estimates again, we obtain $u\in C^\infty(M)$. 
\qed
\begin{remark}
When $u_+$ and $u_-$ are smooth functions on $M$, the proof of Proposition \ref{prop4.1.} is similar to that in \cite[Theorem 2.10]{[Ni]}. In fact, in this case, the sequence $\{v_-^m\}_{m\geq 1}$ defined by \eqref{4.8'} belongs to $C^\infty(\overline{\Omega}_k)$. Hence, we can apply the classical maximum principle (instead of the weak maximum principle) to \eqref{4.9'} to obtain the monotone sequence \eqref{4.13'''}. The remainder of the proof is identical to the one above.
\end{remark}

Applying Proposition \ref{prop4.1.} to the prescribed Chern scalar curvature equation \eqref{2.7...}, we can establish the following existence result.

\begin{theorem}\label{theorem4.2}
Let $(M^n, \omega)$ be a complete noncompact Hermitian manifold with its Chern scalar curvature $S^{Ch}(\omega)$ satisfying
\begin{align}\label{4.6}
-C_1\leq S^{Ch}(\omega)\leq -C_2,
\end{align}
where $C_1\geq C_2>0$ are two constants. Suppose that $K$ is a smooth function on $M$ with 
\begin{equation}\label{4.7}
\begin{aligned}
&K(x)\leq 0, \quad \forall\ x\in M,\\
-a^2\leq & K(x)\leq -b^2, \quad\forall\ x\in M\setminus D,
\end{aligned}
\end{equation}
where $D$ is a compact subset of $M$, $a\geq b>0$ are two constants. Then there exists a complete Hermitian metric $\tilde{\omega}$ which is conformal and uniformly equivalent to $\omega$, such that its Chern scalar curvature $S^{Ch}(\tilde{\omega})=K$.
\end{theorem}
\proof We construct an upper solution $u_+$ and a lower solution $u_-$ of
\begin{align}\label{4.8..}
-\Delta^{Ch}_\omega u+k=Ke^{\frac{2}{n}u},
\end{align}
which satisfy $u_+\geq u_-$ and are bounded on $M$, where $k=S^{Ch}(\omega)$. Since $D$ is compact and $K\in C^\infty(M)$ satisfies \eqref{4.7}, we have $0>\min_{x\in D} K(x)>-\infty$. Set $\tilde{a}^2=\max\{a^2, -\min_{x\in D} K(x)\}$. Then we have 
\begin{align}\label{4.8.}
K(x)\geq -\tilde{a}^2,\quad \forall\ x\in M.
\end{align}
Let $u_-=\frac{n}{2}\log\frac{C_2}{\tilde{a}^2}$. From \eqref{4.6} and \eqref{4.8.}, we obtain
\begin{equation}
\begin{aligned}
-\Delta^{Ch}_\omega u_-+k-Ke^{\frac{2}{n}u_-}&=k-Ke^{\frac{2}{n}\cdot \frac{n}{2}\log\frac{C_2}{\tilde{a}^2}}\\
&\leq -C_2+\tilde{a}^2\frac{C_2}{\tilde{a}^2}\\
&=0.
\end{aligned}
\end{equation}
Next we construct an upper solution 
\begin{align}\label{3.25--}
u_+=\varphi u_0+C, 
\end{align}
where $\varphi$, $u_0$, $C$ will be defined later. Denote
\begin{align}
D_1=\{x\in M: -\frac{1}{2}b^2\leq K\leq 0\}.
\end{align}
Since $K\in C^\infty(M)$ and it satisfies \eqref{4.7}, we obtain that $D_1$ is compact and lies in the interior of $D$. Let $\varphi\in C^\infty_0(M)$ be a cutoff function such that $0\leq \varphi\leq 1$, $\varphi|_{B_1}\equiv 1$ and $\text{supp} \varphi \subset B_{2}\subset D$, where $B_1, B_2$ are open subsets of $M$ satisfying
\begin{align}
D_1\subset B_1\subset \text{supp} \varphi \subset B_2\subset D.
\end{align}
Let $D_2$ be a relatively compact open subset of $M$ with smooth boundary such that $D\subset D_2\subset\subset M$. Let $u_1\in C^\infty(\overline{D_2})$ be the solution of 
\begin{equation}\label{4.11.}
\begin{cases}
-\Delta^{Ch}_\omega u=C_1\quad \text{in}\ D_2\\
u|_{\partial D_2}=0.
\end{cases}
\end{equation}
Here the existence of \eqref{4.11.} is due to Theorem 1.3, Proposition 1.9 and Proposition 2.3 in Chapter 5 of \cite{[Tay]}. Indeed, $-\Delta^{Ch}_\omega =-\Delta_d+\theta^\#$, where $\theta^\#$ is the (real) dual vector field of the Lee form $\theta$ with respect to $\omega$, defined by $\theta^\#u=(du, \theta)_\omega$. By the regularity results (cf. \cite[Chapter 5, Theorem 1.3]{[Tay]}) and the strong maximum principle (cf. \cite[Chapter 5, Proposition 2.3]{[Tay]}) of $-\Delta_d+\theta^\#$, the map $-\Delta^{Ch}_\omega: W^{1,2}_0(D_2)\rightarrow W^{-1,2}(D_2)$ is injective. Consequently, it is also surjective by \cite[Chapter 5, Proposition 1.9]{[Tay]}. Since $D_2$ is a relatively compact open subset with smooth boundary and $C_1$ is a constant, it follows from \cite[Chapter 5, Theorem 1.3]{[Tay]} that the solution $u_1$ of \eqref{4.11.} belongs to $C^\infty(\overline{D_2})$. Now we set
\begin{equation}\label{4.14'}
u_0(x)=
\begin{cases}
u_1,\quad &x\in D_2,\\
0,\quad &x\in M\setminus D_2.
\end{cases}
\end{equation}
It is clear that $u_0\in C^0(M)\cap W^{1,2}_{loc}(M)$. Moreover, $u_0\geq 0$. Indeed, by \eqref{4.11.}, we see that $\Delta^{Ch}_\omega u_1=-C_1<0$ in $D_2$. According to the maximum principle, we have $\min_{\overline{D_2}}u_1=\min_{\partial D_2}u_1=0$, and thus $u_1\geq 0$ in $D_2$. Therefore, $u_0\geq 0$ in $M$. For the constant $C$ in \eqref{3.25--}, we choose
\begin{align}\label{3.30-}
C>\max\left\{\frac{n}{2}\log \frac{2(C'+C_1)}{b^2}, \frac{n}{2}\log\frac{C_2}{\tilde{a}^2}\right\},
\end{align}
where $C'=\max_{\overline{B_{2}}}|u_1\Delta^{Ch}_\omega \varphi+2\nabla \varphi \cdot\nabla u_1|$. Hence, it is obvious that $u_+=\varphi u_0+C$ is a smooth function on $M$ satisfying $u_+\geq C\geq u_-$ on $M$. 

We claim that $u_+=\varphi u_0+C$ is an upper solution of \eqref{4.8..}. Indeed, if $x\in B_{1}$, then using $\varphi|_{B_1}\equiv 1$, \eqref{4.6}, \eqref{4.7} and \eqref{4.11.}, we have for any  $x\in B_1$
\begin{equation}
\begin{aligned}
-\Delta^{Ch}_\omega u_++k-Ke^{\frac{2}{n}u_+}&=-\Delta^{Ch}_\omega u_0+k-Ke^{\frac{2}{n}(u_0+C)}\\
&=-\Delta^{Ch}_\omega u_1+k-Ke^{\frac{2}{n}( u_0+C)}\\
&\geq C_1-C_1\\
&=0.
\end{aligned}
\end{equation}
If $x\in B_2\setminus D_1\subset \{x\in M: K(x)< -\frac{1}{2}b^2\}$, then using \eqref{4.6}, \eqref{4.11.}, \eqref{3.30-}, $0\leq \varphi\leq 1$ and $u_1(x)\geq 0$, we have for any $x\in B_2\setminus D_1$
\begin{equation}
\begin{aligned}
-\Delta^{Ch}_\omega u_++k-Ke^{\frac{2}{n}u_+}&=-u_1\Delta^{Ch}_\omega \varphi-2\nabla \varphi \cdot\nabla u_1-\varphi\Delta^{Ch}_\omega u_1+k-Ke^{\frac{2}{n}u_+}\\
&=-u_1\Delta^{Ch}_\omega \varphi-2\nabla \varphi \cdot\nabla u_1+C_1\varphi+k-Ke^{\frac{2}{n}u_+}\\
&\geq -C'-C_1+\frac{1}{2}b^2e^{\frac{2}{n}\varphi u_1}e^{\frac{2}{n}C}\\
&\geq -C'-C_1+\frac{1}{2}b^2e^{\frac{2}{n}C}\\
&\geq 0,
\end{aligned}
\end{equation}
where $C'=\max_{\overline{B_2}}|u_1\Delta^{Ch}_\omega \varphi+2\nabla \varphi \cdot\nabla u_1|$. If $x\in M\setminus \text{supp} \varphi$, then $u_+(x)=C$, $K(x)<-\frac{1}{2}b^2$, and thus for any  $x\in M\setminus \text{supp} \varphi$
\begin{equation}
\begin{aligned}
-\Delta^{Ch}_\omega u_++k-Ke^{\frac{2}{n}u_+}> -C_1+\frac{1}{2}b^2e^{\frac{2}{n}C}\geq 0,
\end{aligned}
\end{equation}
where we have used \eqref{4.6} and \eqref{3.30-}. Therefore, $u_+$ is a smooth function on $M$ with $-\Delta^{Ch}_\omega u_++k-Ke^{\frac{2}{n}u_+}\geq 0$ in $M$. Moreover,  by \eqref{3.30-}, we have $u_+=\varphi u_0+C\geq C\geq u_-$, and $\sup_M u_+\leq \|u_1\|_{C^0(\overline{D_2})}+C<+\infty$. In terms of Proposition \ref{prop4.1.}, we deduce that there exists a smooth bounded solution of \eqref{4.8..}.
\qed

For the uniqueness, it will be proved by using the generalized maximum principle on complete Hermitian manifolds. First, we introduce the following Laplacian comparison theorem.
\begin{lemma}\label{lemma4.4}
Let $(M^n, \omega)$ be a complete Hermitian manifold with 
\begin{align}\label{4.16}
Ric^{(2)}(X, \overline{X})\geq-C_1(1+r(x))^\alpha,
\end{align}
and the torsion 
\begin{align}\label{4.17.}
\|T^{Ch}(X,Y)\|\leq C_2(1+r(x))^\beta,
\end{align}
where $X, Y\in T^{1,0}_{x}M$ satisfy $\|X\|=\|Y\|=1$, $r(x)$ is the Riemannian distance between a fixed point $x_0$ and $x$ in $(M, \omega)$, and $C_1, C_2, \alpha, \beta$ are positive constants. Let 
\begin{align}
&G(t)=\frac{C_1}{4n}(1+t)^\alpha+\frac{nC_2^2}{2}(1+t)^{2\beta},\\
&h(t)=\frac{1}{\sqrt{G(0)}}\left(e^{\int^t_0G^{\frac{1}{2}}(s)ds}-1\right),\quad t\geq 0.
\end{align}
Then we have
\begin{align}
\Delta^{Ch}_\omega r(x)\leq  4n \frac{h'(r(x))}{h(r(x))}=4n\frac{e^{\int^{r(x)}_0G^{\frac{1}{2}}(s)ds}}{e^{\int^{r(x)}_0G^{\frac{1}{2}}(s)ds}-1}G^{\frac{1}{2}}(r(x))
\end{align}
outside the cut locus of $x_0$.
\end{lemma}
\proof The proof strategy of this theorem is taken from \cite[Lemma 2.1]{[PRS]}, but the details are different in the Hermitian setting. Hence, for the reader's convenience, we provide a full argument. Fix any $x\in M\setminus (\{x_0\}\cup Cut(x_0))$ and let $\gamma(t): [0, l]\rightarrow M$ be a minimal geodesic with respect to the Levi-Civita connection such that $\gamma(0)=x_0$, $\gamma(l)=x$ and $\|\gamma'(0)\|=1$. Let $\{e_i\}_{i=1}^n$ be a unitary frame field of $T^{1,0}M$ on a small neighborhood of $\gamma$. Since $\|\nabla r\|=1$ holds outside the cut locus of $x_0$, using \cite[(3.5)]{[Yu0]}, we have
\begin{equation}\label{4.34''}
\begin{aligned}
0&=\frac{1}{2}\Delta^{Ch}_\omega\|\nabla r\|^2\\
&=\sum_k\left(2\sum_ir_ir_{\bar{i}}\right)_{k\bar{k}}\\
&=2\sum_{k,i}(r_{k\bar{k}i}r_{\bar{i}}+r_{k\bar{k}\bar{i}}r_{i})+2\sum_{i,k}|r_{k\bar{i}}|^2+2\sum_{i,j,k}\left(r_{\bar{i}}\tau^j_{ik}r_{j\bar{k}}+r_{k\bar{j}}\tau^{\bar{j}}_{\bar{i}\bar{k}}r_i \right)\\
&\quad+2\sum_{i,k}|r_{ik}|^2+2\sum_{i,j,k} R_{i\bar{j}k\bar{k}}r_{j}r_{\bar{i}},
\end{aligned}
\end{equation}
where $T^{Ch}(e_i, e_k)=\sum_j\tau^j_{ik}e_j$. By \eqref{4.16}, \eqref{4.17.}, $\|\nabla r\|^2=2\sum_i|r_i|^2=1$, and using the Cauchy-Schwarz and Young inequalities, we have
\begin{equation}\label{3.40+}
\begin{aligned}
\sum_{i,j,k} R_{i\bar{j}k\bar{k}}r_{j}r_{\bar{i}}&= Ric^{(2)}(\sum_i r_{\bar{i}}e_i, \sum_j r_{j}e_{\bar{j}}) \geq\frac{1}{2}C_1(1+r(x))^\alpha,\\
|\sum_{i,j,k}r_{\bar{i}}\tau^j_{ik}r_{j\bar{k}}|&=|\sum_{j,k} r_{j\bar{k}}T^{Ch}(\sum_i r_{\bar{i}}e_i, e_k) |\\
&\leq \frac{1}{\sqrt{2}}C_2(1+r(x))^\beta\sum_{j,k}|r_{j\bar{k}}|\\
&\leq \frac{1}{\sqrt{2}}C_2(1+r(x))^\beta n(\sum_{j,k}|r_{j\bar{k}}|^2)^{\frac{1}{2}}\\
&\leq \frac{1}{2}C_2^2(1+r(x))^{2\beta} n^2+\frac{1}{4}\sum_{j,k}|r_{j\bar{k}}|^2,\\
\sum_{i,k}|r_{k\bar{i}}|^2&\geq \sum_{k}|r_{k\bar{k}}|^2\geq \frac{1}{n}\left(\sum_{k}r_{k\bar{k}}\right)^2=\frac{1}{4n}(\Delta^{Ch}_\omega r)^2.
\end{aligned}
\end{equation}
Similarly,
\begin{align}\label{3.41+}
|\sum_{i,j,k}r_{k\bar{j}}\tau^{\bar{j}}_{\bar{i}\bar{k}}r_i|\leq  \frac{1}{2}C_2^2(1+r(x))^{2\beta} n^2+\frac{1}{4}\sum_{j,k}|r_{k\bar{j}}|^2.
\end{align}
Hence, it follows from \eqref{4.34''}, \eqref{3.40+}, \eqref{3.41+} that
\begin{align}\label{4.21.}
0\geq \langle \nabla\Delta^{Ch}_\omega r, \nabla r\rangle+\frac{1}{4n}(\Delta^{Ch}_\omega r)^2-4nG(r),
\end{align}
where 
\begin{align}
G(r)=\frac{C_1}{4n}(1+r)^\alpha+\frac{nC_2^2}{2}(1+r)^{2\beta}.
\end{align}
Set 
\begin{align}
f(t)=\Delta^{Ch}_\omega r(\gamma(t)),\quad t\in (0,l].
\end{align}
Then \eqref{4.21.} can be written as
\begin{align}\label{4.24}
0\geq \frac{1}{4n} f^2(t)+f'(t)-4nG(t).
\end{align}
Recall a well-known result in Riemannian geometry:
\begin{align}
\Delta_d r=\frac{2n-1}{r}+o(1),\ \text{as}\ r\rightarrow 0+.
\end{align}
Since $\|dr\|\equiv 1$ on $M\setminus \left(\{x_0\}\cup \text{Cut}(x_0)\right)$, we have
\begin{align}\label{4.25.}
f(t)= \Delta_d r(\gamma(t))-(dr,\theta)_{\omega}(\gamma(t))=\frac{2n-1}{t}+O(1),\ \text{as}\ t\rightarrow 0+,
\end{align}
which implies the function 
\begin{equation}
g(t)=t^ae^{\int_0^t(\frac{f(s)}{4n}-\frac{a}{s})ds},\quad a=\frac{2n-1}{4n}
\end{equation}
 is well-defined on $[0, l]$. A simple computation yields
 \begin{align}\label{4.27}
 g(0)=0, \quad g'(t)=\frac{f}{4n}g \ \  (t\in (0, l]),
 \end{align}
 and by \eqref{4.24} and \eqref{4.27}, 
 \begin{align}\label{4.28}
 g''(t)\leq G(t)g(t).
 \end{align}
Let
\begin{align}\label{4.29}
h(t)=\frac{1}{\sqrt{G(0)}}\{e^{\int^t_0G^{\frac{1}{2}}(s)ds}-1\},\quad t\geq 0.
\end{align}
Then 
\begin{align}\label{4.30}
h''(t)-G(t)h(t)\geq 0,\quad h(0)=0, \quad h'(0)=1.
\end{align}
Set
\begin{align}
\Phi(t)=(gh'-g'h)(t).
\end{align}
Using \eqref{4.28}, \eqref{4.30} and $h\geq 0$, we have
\begin{align}
\Phi'(t)=gh''-g''h\geq 0.
\end{align}
Hence,
\begin{align}\label{4.33}
\Phi(t)\geq \lim_{t\rightarrow 0+}\Phi(t)=0.
\end{align}
Indeed, by \eqref{4.25.}-\eqref{4.27} and \eqref{4.29}, we deduce that
\begin{equation}
\begin{aligned}
\lim_{t\rightarrow 0+}g'(t)h(t)&=\frac{1}{4n\sqrt{G(0)}}\lim_{t\rightarrow 0+}\left(tf(t)e^{\int_0^t(\frac{f(s)}{4n}-\frac{a}{s})ds}\frac{e^{\int^t_0G^{\frac{1}{2}}(s)ds}-1}{t^{1-a}}\right)\\
&=0,
\end{aligned}
\end{equation}
and thus
\begin{align}
\lim_{t\rightarrow 0+}\Phi(t)=(g(0)h'(0)-\lim_{t\rightarrow 0+}g'(t)h(t))=0.
\end{align}
According to \eqref{4.27} and \eqref{4.33}, we get
\begin{align}
f(t)=4n \frac{g'(t)}{g(t)}\leq 4n \frac{h'(t)}{h(t)}.
\end{align}
\qed

In particular, if we choose $\alpha=2$ and $\beta=1$ in Lemma \ref{lemma4.4}, then when $x$ is a smooth point of $r$ and is outside some compact subset of $M$, we have
\begin{equation}
\begin{aligned}
&\Delta^{Ch}_\omega r^2(x)=2r\Delta^{Ch}_\omega r+2\|\nabla r\|^2\leq C_4r(x)G^{\frac{1}{2}}(r(x)), \quad G(r)=C_3(1+r)^2\\
 &\|\nabla r^2\|(x)=2r(x)
\end{aligned}
\end{equation}
for some constants $C_3, C_4>0$. It is easy to check that $G$ satisfies the condition \cite[(1.17)]{[PRS]}. Hence, following the proof of \cite[Theorem 1.9, Remark 1.11]{[PRS]}, we can easily obtain the following generalized maximum principle. In fact, in the proof of \cite[Theorem 1.9]{[PRS]}, they use the function $G$ to construct an auxiliary function that attains a maximum on $M$ and then apply the classical maximum principle to obtain their results. In our case, the classical maximum principle still holds for the Chern Laplacian because the first-order term of the Chern Laplacian vanishes at the maximum point. Therefore, following the same argument, we can also prove the following generalized maximum principle on Hermitian manifolds.

\begin{proposition}\label{prop4.1}
Let $(M^n, \omega)$ be a complete Hermitian manifold with 
\begin{align}
Ric^{(2)}(X, \overline{X})\geq-C_1(1+r(x))^2,
\end{align}
and the torsion 
\begin{align}
\|T^{Ch}(X,Y)\|\leq C_2(1+r(x)),
\end{align}
where $X, Y\in T^{1,0}_{x}M$ satisfy $\|X\|=\|Y\|=1$, $r(x)$ is the Riemannian distance between a fixed point $x_0$ and $x$ in $(M, \omega)$, and $C_1, C_2$ are positive constants. Then for any $C^2$ function $f:M\rightarrow \mathbb{R}$ with $\sup_M f<+\infty$, there exists a sequence $\{x_m\}\subset M$ such that
\begin{align}
\lim_{m\rightarrow \infty}f(x_m)=\sup_Mf,\quad \lim_{m\rightarrow \infty}\|\nabla f\|(x_m)=0,\quad \limsup_{m\rightarrow \infty}\Delta^{Ch}_\omega f(x_m)\leq 0.
\end{align}
\end{proposition}
\begin{remark}
It is well-known that the generalized maximum principle due to Omori \cite{[Omo]} and Yau \cite{[Yau]} is an important analytical tool in differential geometry and has lots of applications (cf. \cite{[Yau2]}, \cite{[PRS]}, \cite{[Sun]}, etc.).
\end{remark}

Making use of the above generalized maximum principle, we can prove the following uniqueness result.

\begin{theorem}\label{theorem4.6}
Let $(M^n, \omega)$ be a complete noncompact Hermitian manifold with 
\begin{align}
Ric^{(2)}(X, \overline{X})\geq-C_1(1+r(x))^2,
\end{align}
and the torsion 
\begin{align}
\|T^{Ch}(X,Y)\|\leq C_2(1+r(x)),
\end{align}
where $X, Y\in T^{1,0}_{x}M$ satisfy $\|X\|=\|Y\|=1$, $r(x)$ is the Riemannian distance between a fixed point $x_0$ and $x$ in $(M, \omega)$, and $C_1, C_2$ are positive constants. Suppose that $K$ is a smooth function on $M$ with 
\begin{equation}\label{4.57'''}
\begin{aligned}
&K(x)\leq 0, \quad \forall\ x\in M,\\
&K(x)\leq -b^2, \quad\forall\ x\in M\setminus D,
\end{aligned}
\end{equation}
where $D$ is a compact subset of $M$, $b>0$ is a constant. Then there exists at most one bounded smooth solution of 
\begin{align}\label{4.8}
-\Delta^{Ch}_\omega u+k=Ke^{\frac{2}{n}u},
\end{align}
where $k\in C^\infty(M)$ is an arbitrary function on $M$.
\end{theorem}
\proof Let $u_1, u_2$ be two bounded smooth solutions of \eqref{4.8} and set $w=u_1-u_2$. Then $w$ is bounded and satisfies
\begin{align}\label{4.59'}
-\Delta^{Ch}_\omega w=K\left(e^{\frac{2}{n}u_1}-e^{\frac{2}{n}u_2}\right).
\end{align}
According to Proposition \ref{prop4.1}, there exists a sequence $\{x_m\}\subset M$ such that
\begin{align}\label{4.10}
\lim_{m\rightarrow \infty}w(x_m)=\sup_Mw,\quad \lim_{m\rightarrow \infty}\|\nabla w\|(x_m)=0,\quad \limsup_{m\rightarrow \infty}\Delta^{Ch}_\omega w(x_m)\leq 0.
\end{align}
\textbf{Case 1:} There exists a subsequence $\{x_{m_k}\}$ and $x_0\in M$ such that $x_{m_k}\rightarrow x_0$ as $k\rightarrow \infty$.

By \eqref{4.10}, we have $w(x_0)=\sup_Mw$ and $\|\nabla w\|(x_0)=0$. Now we prove that $w\leq 0$ in $M$. Suppose that $w(x_0)=\sup_Mw>0$, then $w>0$ in a neighborhood $U$ of $x_0$. It follows from \eqref{4.57'''} and \eqref{4.59'} that $\Delta^{Ch}_\omega w\geq 0$ in $U$. Since $x_0$ is a maximum point of $w$ in $U$,  by the strong maximum principle we obtain $w$ is constant in $U$. Enlarging $U$, we get that $w$ is a positive constant on $M$, and thus $K\equiv 0$ because of \eqref{4.59'}, which leads to a contradiction. Therefore, $w\leq 0$ in $M$. \\
\textbf{Case 2:} $\{x_{m}\}$ does not have an accumulation point in $M$.

Assume that $\lim_{m\rightarrow \infty}w(x_m)=\sup_Mw>0$. Then for any small $\epsilon>0$ there exists a constant $N>0$ such that $w(x_m)>\epsilon>0$ for all $m\geq N$. Hence,
\begin{equation}\label{4.11}
\begin{aligned}
e^{\frac{2}{n}u_1(x_m)}-e^{\frac{2}{n}u_2(x_m)}&=e^{\frac{2}{n}u_2(x_m)}\left(e^{\frac{2}{n}w(x_m)}-1 \right)\\
&>e^{\frac{2}{n}\inf_M u_2}\left(e^{\frac{2}{n}\epsilon}-1 \right)\\
&>0
\end{aligned}
\end{equation}
for any $m\geq N$, where we have used the assumption $u_2$ is a bounded solution on $M$. Since $\{x_m\}$ does not have an accumulation point in $M$, we obtain $r(x_m)\rightarrow +\infty$ as $m\rightarrow +\infty$.
Hence, from \eqref{4.57'''} and \eqref{4.11} it follows that for sufficiently large $m$,
\begin{equation}
\begin{aligned}
\Delta^{Ch}_\omega w(x_m)&=-K\left(e^{\frac{2}{n}u_1(x_m)}-e^{\frac{2}{n}u_2(x_m)}\right)\\
&>b^2e^{\frac{2}{n}\inf_M u_2}\left(e^{\frac{2}{n}\epsilon}-1 \right)>0,
\end{aligned}
\end{equation}
which leads to a contradiction with $\limsup_{m\rightarrow \infty}\Delta^{Ch}_\omega w(x_m)\leq 0$. Therefore, $w\leq \sup_M w\leq 0$.

In conclusion, we get $w\leq 0$ in $M$. Similarly, we also have $w\geq 0$ in $M$. Thus $w\equiv 0$ on $M$, and consequently $u_1=u_2$ on $M$.
\qed

By combining Theorem \ref{theorem4.2} and Theorem \ref{theorem4.6}, we obtain

\begin{theorem}\label{theorem4.7}
Let $(M^n, \omega)$ be a complete noncompact Hermitian manifold with 
\begin{align}\label{3.69+}
Ric^{(2)}(X, \overline{X})\geq-C_1(1+r(x))^2,
\end{align}
and the torsion 
\begin{align}\label{3.70+}
\|T^{Ch}(X,Y)\|\leq C_2(1+r(x)),
\end{align}
where $X, Y\in T^{1,0}_{x}M$ satisfy $\|X\|=\|Y\|=1$, $r(x)$ is the Riemannian distance between a fixed point $x_0$ and $x$ in $(M, \omega)$, $C_1, C_2$ are positive constants. In addition, we assume that its Chern scalar curvature $S^{Ch}(\omega)$ satisfies
\begin{align}
-C_3\leq S^{Ch}(\omega)\leq -C_4,
\end{align}
where $C_3, C_4$ are positive constants. Suppose that $K$ is a smooth function on $M$ with 
\begin{equation}\label{4.}
\begin{aligned}
&K(x)\leq 0, \quad \forall\ x\in M,\\
-a^2\leq & K(x)\leq -b^2, \quad\forall\ x\in M\setminus D,
\end{aligned}
\end{equation}
where $D$ is a compact subset of $M$, $a\geq b>0$ are constants. Then there exists a unique complete Hermitian metric $\tilde{\omega}$ which is conformal and uniformly equivalent to $\omega$, such that its Chern scalar curvature $S^{Ch}(\tilde{\omega})=K$.
\end{theorem}

\proof[Proof of Theorem \ref{theorem1.3}] The theorem follows from Theorem \ref{theorem4.2} and Theorem \ref{theorem4.7}. 
\qed

Making use of Theorem \ref{theorem1.3} and Theorem \ref{theorem4.6}, we can also prove the following existence and uniqueness results.

\begin{theorem}[=Theorem \ref{theorem1.6}]\label{theorem3.8-}
Let $(M^n, \omega)$ be a complete noncompact Hermitian manifold with its Chern scalar curvature $S^{Ch}(\omega)$ satisfying
\begin{equation}\label{3.80...}
\begin{aligned}
&S^{Ch}(\omega)(x)\leq 0, \quad\forall\ x\in M\\
-a^2\leq &S^{Ch}(\omega)(x)\leq -b^2, \quad\forall\ x\in M\setminus D,
\end{aligned}
\end{equation}
where $D$ is a compact subset of $M$, $a\geq b>0$ are two constants. Then there exists a complete Hermitian metric $\tilde{\omega}$ which is conformal and uniformly equivalent to $\omega$, such that its Chern scalar curvature $S^{Ch}(\tilde{\omega})\equiv-1$. If, in addition, we assume that $(M, \omega)$ satisfies \eqref{3.69+} and \eqref{3.70+}, then $\tilde{\omega}$ is unique among all conformal metrics that are uniformly equivalent to $\omega$.
\end{theorem}

\proof Let $\phi\in C^\infty_0(M)$ satisfy $\Delta^{Ch}_\omega\phi=\delta$ in $D$ and $\Delta^{Ch}_\omega\phi\geq-\frac{1}{2}b^2$ in $M$, where $\delta>0$ is a constant.  Indeed, let $D_1, D_2$ be two relatively compact open subsets of $M$ with smooth boundaries such that $D\subset D_1\subset\subset D_2$. Consider the function $\phi_0\in C^\infty(\overline{D_2})$ satisfying
\begin{align}
\begin{cases}
&\Delta^{Ch}_\omega \phi_0=\delta_1 \quad\text{in}\ D_2,\\
& \phi_0|_{\partial D_2}=\delta_2,
\end{cases}
\end{align}
where $\delta_1, \delta_2>0$ are constants. Extend $\phi_0|_{D_1}$ smoothly to a $C^\infty$ function on $M$ with compact support, and denote this extension by $\phi_1$. Then $\phi_1\in C^\infty_0(M)$ and it satisfies $\Delta^{Ch}_\omega \phi_1=\delta_1$ in $D$ and $\Delta^{Ch}_\omega \phi_1\geq-C$ in $M$ for some constant $C>0$. Set $\phi=\frac{b^2}{2C}\phi_1\in  C^\infty_0(M)$. Then $\Delta^{Ch}_\omega \phi=\delta$ in $D$ and $\Delta^{Ch}_\omega \phi\geq -\frac{b^2}{2}$ in $M$, where $\delta=\frac{b^2}{2C}\delta_1>0$. 

Let $\omega_1=e^{\frac{2}{n}\phi}\omega$. Then we have
\begin{equation}\label{3.75--}
\begin{aligned}
S^{Ch}(\omega_1)=e^{-\frac{2}{n}\phi}(-\Delta^{Ch}_\omega \phi+S^{Ch}(\omega))&\leq-C_2\quad \text{in}\ M,
\end{aligned}
\end{equation}
where we have used \eqref{3.80...} and $S^{Ch}(\omega)\leq 0$ in $M$, and $C_2=e^{-\frac{2}{n}\|\phi\|_{C^0(M)}}\min\{\delta, \frac{1}{2}b^2\}$. Moreover, since $D$ is compact and $S^{Ch}(\omega)\in C^\infty(M)$ satisfies \eqref{3.80...}, we have $0>\min_{x\in D} S^{Ch}(\omega)(x)>-\infty$. Set $\tilde{a}^2=\max\{a^2, -\min_{x\in D} S^{Ch}(\omega)(x)\}$, then we obtain 
\begin{align}\label{4.8.+}
S^{Ch}(\omega)(x)\geq -\tilde{a}^2,\quad \forall\ x\in M.
\end{align}
Hence, 
\begin{equation}\label{3.77--}
\begin{aligned}
S^{Ch}(\omega_1)&=e^{-\frac{2}{n}\phi}(-\Delta^{Ch}_\omega \phi+S^{Ch}(\omega))\\
&\geq -e^{\frac{2}{n}\|\phi\|_{C^0(M)}}(\|\Delta^{Ch}_\omega \phi\|_{C^0(M)}+\tilde{a}^2).
\end{aligned}
\end{equation}
Combining \eqref{3.75--} and \eqref{3.77--} yields
\begin{align}
-C_1\leq S^{Ch}(\omega_1)\leq -C_2\quad \text{in}\ M,
\end{align}
where $C_1=e^{\frac{2}{n}\|\phi\|_{C^0(M)}}(\|\Delta^{Ch}_\omega \phi\|_{C^0(M)}+\tilde{a}^2)$. Since $K\equiv-1$ satisfies the condition \eqref{4.7} and $S^{Ch}(\omega_1)$ satisfies \eqref{4.6}, by Theorem \ref{theorem4.2} there exists a complete Hermitian metric $\tilde{\omega}$ that is conformal and uniformly equivalent to $\omega_1$ and satisfies $S^{Ch}(\tilde{\omega})\equiv-1$. From $\omega_1=e^{\frac{2}{n}\phi}\omega$ and $\phi\in C^\infty_0(M)$, it follows that the complete Hermitian metric $\tilde{\omega}$ is also conformal and uniformly equivalent to $\omega$. Therefore, we have completed the proof of the existence result in Theorem \ref{theorem3.8-}.

For the uniqueness, since $K=S^{Ch}(\tilde{\omega})\equiv-1$ satisfies the condition \eqref{4.57'''}, Theorem \ref{theorem4.6} implies that such $\tilde{\omega}$ is unique among all conformal metrics that are uniformly equivalent to $\omega$.
\qed

Finally, we study the case where $S^{Ch}(\omega)$ is not bounded below by $-a^2$ outside some compact subset $D$ of $M$. Before presenting our theorem, we prove the following lemma using the classical maximum principle.

\begin{lemma}\label{lemma4.9}
Let $(M^n, \omega)$ be a complete Hermitian manifold. Suppose that $u\in C^\infty(M)$ satisfies
\begin{align}\label{4.61.}
-\Delta^{Ch}_{\omega} u+S\leq-e^{\frac{2}{n}u},
\end{align}
where $S\in C^\infty(M)$. Then for any compact subset $K\subset M$, there exists a constant $C=C(n, K, S)$ such that
\begin{align}
\sup_K u<C.
\end{align}
\end{lemma}
\proof Set 
\begin{align}
v=e^{\frac{2}{n}u}.
\end{align}
By \eqref{4.61.}, we obtain
\begin{align}\label{4.64..}
-\Delta^{Ch}_\omega v+\frac{|\nabla v|^2}{v}+\frac{2}{n}Sv\leq-\frac{2}{n}v^2.
\end{align}
We choose a cutoff function $\varphi\in C^\infty_0(M)$ with $0\leq \varphi\leq 1$ on $M$ and $\varphi\equiv 1$ on $K$. Set
\begin{align}
\tilde{v}=\varphi^2v=\varphi^2e^{\frac{2}{n}u}.
\end{align}
Then there exists a point $x_0\in \text{supp} \varphi$ such that $\tilde{v}(x_0)=\max_M \tilde{v}>0$. By the maximum principle, we have
\begin{align}\label{4.66..}
0=\nabla \tilde{v}(x_0)=\nabla(\varphi^2v)=2\varphi v \nabla \varphi +\varphi^2\nabla v,
\end{align}
and
\begin{equation}\label{4.67..}
\begin{aligned}
0&\geq \Delta^{Ch}_\omega \tilde{v}(x_0)\\
&=\Delta^{Ch}_\omega (\varphi^2v)\\
&=2|\nabla \varphi|^2v+2\varphi v\Delta^{Ch}_\omega\varphi +4\varphi\nabla\varphi \cdot \nabla v+\varphi^2\Delta^{Ch}_\omega v.
\end{aligned}
\end{equation}
From \eqref{4.66..}, we deduce that at $x_0$
\begin{align}\label{4.68...}
\nabla v=-2\varphi^{-1}v \nabla \varphi.
\end{align}
Substituting \eqref{4.64..} and \eqref{4.68...} into \eqref{4.67..} yields, at $x_0$,
\begin{equation}
\begin{aligned}
0&\geq2|\nabla \varphi|^2v+2\varphi v\Delta^{Ch}_\omega\varphi +4\varphi\nabla\varphi \cdot \nabla v+\varphi^2\Delta^{Ch}_\omega v\\
&\geq2|\nabla \varphi|^2v+2\varphi v\Delta^{Ch}_\omega\varphi +4\varphi\nabla\varphi \cdot (-2\varphi^{-1}v \nabla \varphi)\\
&\quad\quad+\varphi^2\left( \frac{|-2\varphi^{-1}v \nabla \varphi|^2}{v}+\frac{2}{n}Sv+\frac{2}{n}v^2\right)\\
&=-2|\nabla \varphi|^2v+2\varphi v\Delta^{Ch}_\omega\varphi +\varphi^2\frac{2}{n}Sv+\frac{2}{n}\varphi^2v^2,
\end{aligned}
\end{equation}
which implies that
\begin{equation}
\begin{aligned}
\max_M \tilde{v}=\tilde{v}(x_0)&=(\varphi^2v)(x_0)\\
&\leq -\varphi^2(x_0)S(x_0)-n\varphi(x_0)\Delta^{Ch}_\omega \varphi(x_0)+n|\nabla \varphi|^2(x_0)\\
&\leq -\varphi^2(x_0)\min_{\text{supp}\varphi}S-n\varphi(x_0)\Delta^{Ch}_\omega \varphi(x_0)+n|\nabla \varphi|^2(x_0)\\
&\leq C
\end{aligned}
\end{equation}
for some constant $C=C(K, \|\varphi\|_{C^2(M)},\min_{\text{supp}\varphi}S)>0$. Hence, for any $x\in K$
\begin{align}
e^{\frac{2}{n}u(x)}= (\varphi^2v)(x)=\tilde{v}(x)\leq \max_M \tilde{v}\leq C,
\end{align}
where we have used $\varphi|_K\equiv 1$. Therefore, we get $\sup_Ku <C$ for some uniform constant $C>0$.
\qed

Following the idea of Aviles--McOwen in \cite{[AM-2]} and using the above lemma, we prove a Hermitian version of their existence result for metrics with constant negative Riemannian scalar curvature, although the corresponding equations differ.

\begin{theorem}[=Theorem \ref{theorem1.7.}]\label{theorem3.10-}
Let $(M^n, \omega)$ be a complete noncompact Hermitian manifold with nonpositive Chern scalar curvature $S^{Ch}(\omega)$ satisfying
\begin{equation}\label{4.61}
\begin{aligned}
S^{Ch}(\omega)(x)\leq -b^2, \quad\forall\ x\in M\setminus D,
\end{aligned}
\end{equation}
where $D$ is a compact subset of $M$, and $b>0$ is a constant. Then there exists a complete Hermitian metric $\tilde{\omega}$ which is conformal to $\omega$ and satisfies $\tilde\omega\geq C\omega$ for some constant $C>0$, such that its Chern scalar curvature $S^{Ch}(\tilde{\omega})\equiv-1$.
\end{theorem}
\proof Firstly, we reduce \eqref{4.61} to the case where \eqref{4.61} holds on all of $M$. As in the proof of Theorem \ref{theorem3.8-}, let $\phi\in C^\infty_0(M)$ satisfy $\Delta^{Ch}_\omega\phi=\delta$ in $D$ and $\Delta^{Ch}_\omega\phi\geq-\frac{1}{2}b^2$ in $M$, where $\delta>0$ is a constant. Let $\omega_1=e^{\frac{2}{n}\phi}\omega$. Then we have
\begin{equation}
\begin{aligned}
S^{Ch}(\omega_1)=e^{-\frac{2}{n}\phi}(-\Delta^{Ch}_\omega \phi+S^{Ch}(\omega))&\leq-\epsilon\quad \text{in}\ M,
\end{aligned}
\end{equation}
where we have used \eqref{4.61} and $S^{Ch}(\omega)\leq 0$ in $M$, and $\epsilon=e^{-\frac{2}{n}\|\phi\|_{C^0(M)}}\min\{\delta, \frac{1}{2}b^2\}$. If $\tilde{\omega}=e^{\frac{2}{n}u}\omega_1$ for some $u\in C^\infty(M)$, then $u$ satisfies
\begin{align}
-\Delta^{Ch}_{\omega_1} u+S^{Ch}(\omega_1)=-e^{\frac{2}{n}u}.
\end{align}
To solve the above equation, we first choose $u_-=a$ where $a<\frac{n}{2}\log \epsilon$ is a constant with 
\begin{align}
-\Delta^{Ch}_{\omega_1} u_-+S^{Ch}(\omega_1)+e^{\frac{2}{n}u_-}=S^{Ch}(\omega_1)+e^{\frac{2}{n}a}\leq -\epsilon+e^{\frac{2}{n}a}<0.
\end{align}

Suppose $M=\cup_{k\geq 1} \Omega_k$, where each $\Omega_k$ is a bounded domain with smooth boundary and $\Omega_k\subset\subset \Omega_{k+1}$. For any $\Omega_k$, there exists a constant $C_k$ dependent on $\Omega_k$ such that $S^{Ch}(\omega_1)>-C_k$. Thus, we can construct an upper solution $u_+=b_k$ in $\Omega_k$, where $b_k$ is a constant satisfying $b_k>\frac{n}{2}\log C_k$, since
\begin{align}
-\Delta^{Ch}_{\omega_1} u_++S^{Ch}(\omega_1)+e^{\frac{2}{n}u_+}=S^{Ch}(\omega_1)+e^{\frac{2}{n}b_k}>-C_k+e^{\frac{2}{n}b_k}>0\ \text{in}\ \Omega_k.
\end{align}
Now we consider the boundary value problem
\begin{equation}\label{4.68}
\begin{cases}
-\Delta^{Ch}_{\omega_1} u+S^{Ch}(\omega_1)=-e^{\frac{2}{n}u}\quad \text{in}\ \Omega_k\\
u=\frac{a+b_k}{2}\quad \text{on}\ \partial \Omega_k
\end{cases}
\end{equation}
By the well-known Monotone Iteration Scheme (see also the proof of Proposition \ref{prop4.1.}), there exists a solution $u_k\in C^\infty(\Omega_k)$ of \eqref{4.68} with $a=u_-\leq u_k\leq u_+=b_k$.

Let us consider the sequence $\{u_k\}_{k\geq 4}$ on $\overline{\Omega}_3$. According to Lemma \ref{lemma4.9}, we have $\sup_{\overline{\Omega}_3}u_k\leq C$ for some constant $C>0$ which is independent of $k$. Together with $u_k\geq u_-=a$, this yields $\|u_k\|_{C^0(\overline{\Omega}_3)}\leq C$. Using the interior elliptic $L^p$-estimates, there is a constant $C>0$ such that
\begin{equation}
\begin{aligned}
\|u_k\|_{W^{2,p}(\Omega_2)}&\leq C\left(\|u_k\|_{L^p(\Omega_3)}+\|S^{Ch}(\omega_1)+e^{\frac{2}{n}u_k}\|_{L^p(\Omega_3)}\right)\\
&\leq  C\left(\|u_k\|_{C^0(\Omega_3)}|\Omega_3|^{\frac{1}{p}}+\|S^{Ch}(\omega_1)\|_{L^p(\Omega_3)}+e^{\frac{2}{n}\|u_k\|_{C^0(\Omega_3)}}|\Omega_3|^{\frac{1}{p}}\right)\\
&\leq C'
\end{aligned}
\end{equation}
for any $p>1$, where $C'$ is a constant independent of $k$. Taking $p=2n+1$ and using the Sobolev embedding theorem $W^{2, 2n+1}(\Omega_2)\subset C^1(\Omega_2)$ yields $\|u_k\|_{C^1(\Omega_2)}\leq C$. According to the Schauder estimates, we have $\|u_k\|_{C^{2, \alpha}(\Omega_1)}\leq C$ for any $\alpha\in (0,1)$. By the compact embedding theorem $C^{2, \alpha}(\Omega_1)\subset\subset C^{2}(\Omega_1)$, it follows that there exists a subsequence $\{u_{k1}\}$ of $\{u_k\}$ such that $u_{k1}$ converges to $u_1$ in $C^2(\Omega_1)$ as $k\rightarrow +\infty$, where $u_1$ solves
\begin{align}
-\Delta^{Ch}_{\omega_1} u+S^{Ch}(\omega_1)=-e^{\frac{2}{n}u}\quad \text{in}\ \Omega_1.
\end{align}
Repeating the above procedure with $\{u_{k1}\}_{k\geq 5}$ on $\Omega_4$, we obtain a subsequence $\{u_{k2}\}$ of $\{u_{k1}\}$ such that  $u_{k2}$ converges to $u_2$ in $C^2(\Omega_2)$ as $k\rightarrow +\infty$, and $u_2$ solves
\begin{equation}
\begin{cases}
-\Delta^{Ch}_{\omega_1} u+S^{Ch}(\omega_1)=-e^{\frac{2}{n}u}\quad \text{in}\ \Omega_2,\\
u|_{\Omega_1}=u_1
\end{cases}
\end{equation}
Inductively, one can obtain a subsequence $\{u_{kj}\}$ which converges to $u_j$ in $C^2(\Omega_j)$ as $k\rightarrow +\infty$. Here $u_j$ is a solution of 
\begin{equation}
\begin{cases}
-\Delta^{Ch}_{\omega_1} u+S^{Ch}(\omega_1)=-e^{\frac{2}{n}u}\quad \text{in}\ \Omega_j,\\
u|_{\Omega_{j-1}}=u_{j-1}
\end{cases}
\end{equation}
for any $j\geq 2$. Define 
\begin{align}
u(x)=\lim_{k\rightarrow +\infty}u_{kk}(x)\quad \ x\in M.
\end{align}
Then $u\in C^2(M)$ and it solves
\begin{equation}
-\Delta^{Ch}_{\omega_1} u+S^{Ch}(\omega_1)=-e^{\frac{2}{n}u}\quad \text{in}\ M.
\end{equation}
By the Schauder estimates, we obtain $u\in C^\infty(M)$. Therefore, there exists a Hermitian metric $\tilde{\omega}$ conformal to $\omega$ such that its Chern scalar curvature $S^{Ch}(\tilde{\omega})=-1$. Moreover, $u\geq u_-=a$. Since $\omega$ is complete and $\tilde{\omega}=e^{\frac{2}{n}u}\omega_1=e^{\frac{2}{n}(u+\phi)}\omega$, where $\phi\in C^\infty_0(M)$ and $u+\phi\geq a-\|\phi\|_{C^0(M)}>-\infty$, it follows that $\tilde{\omega}$ is also complete. 
\qed

Finally, we consider the case that the candidate curvature function is sign-changing. We prove the following theorem:

\begin{theorem}\label{theorem4.8+}
Let $(M^n, \omega)$ be a complete noncompact Hermitian manifold with its Chern scalar curvature $S^{Ch}(\omega)$ satisfying 
\begin{align}\label{4.6-}
-C_1\leq S^{Ch}(\omega)\leq -C_2,
\end{align}
where $C_1, C_2$ are positive constants. Suppose that $K$ is a smooth function on $M$ with 
\begin{equation}\label{3.103}
-a^2\leq  K(x)\leq -b^2, \quad\forall\ x\in M\setminus D,
\end{equation}
where $D$ is a compact subset of $M$, $a\geq b>0$ are constants. There exists a constant $\epsilon>0$ with the following property: if $K(x)\leq \epsilon$ for all $x\in M$, then there exists a complete Hermitian metric $\tilde{\omega}$ that is conformal and uniformly equivalent to $\omega$ and satisfies $S^{Ch}(\tilde{\omega})=K$.
\end{theorem}
\proof We construct an upper solution $u_+$ and a lower solution $u_-$ of
\begin{align}\label{3.104}
-\Delta^{Ch}_\omega u+k=Ke^{\frac{2}{n}u},
\end{align}
which satisfy $u_+\geq u_-$ and are bounded on $M$, where $k=S^{Ch}(\omega)$. Let $u_-$ be as in the proof of Theorem \ref{theorem4.2}. Note that the upper solution $u_+$ constructed in Theorem \ref{theorem4.2} is not applicable to the case of the present theorem. In what follows, we provide a different construction; nevertheless, the upper solution constructed below remains valid for Theorem \ref{theorem4.2}.

Next we construct an upper solution 
\begin{align}\label{3.25---}
u_+=\varphi_1 u_0+C, 
\end{align}
where $\varphi_1$, $u_0$, $C$ are defined later.  Let $\varphi_1\in C^\infty_0(M)$ be a cutoff function such that $0\leq \varphi_1\leq 1$, $\varphi_1|_{B_1}\equiv 1$ and $\text{supp} \varphi_1 \subset B_{2}$, where $B_1, B_2$ are relatively compact open subsets of $M$ satisfying
\begin{align}
D\subset B_1\subset \text{supp} \varphi_1 \subset  B_2\subset\subset M.
\end{align}
Let $D_2$ be a relatively compact open subset of $M$ with smooth boundary such that $B_2\subset D_2\subset\subset M$. Let $u_1\in C^\infty(\overline{D_2})$ be the solution of 
\begin{equation}\label{4.11.-}
\begin{cases}
-\Delta^{Ch}_\omega u=2C_1\quad \text{in}\ D_2\\
u|_{\partial D_2}=0.
\end{cases}
\end{equation}
Now we set
\begin{equation}\label{4.14'-}
u_0(x)=
\begin{cases}
u_1,\quad &x\in D_2,\\
0,\quad &x\in M\setminus D_2.
\end{cases}
\end{equation}
As in the proof of Theorem \ref{theorem4.2}, $u_0\in C^0(M)\cap W^{1,2}_{loc}(M)$ and $u_0\geq 0$. For the constant $C$ in \eqref{3.25---}, we choose
\begin{align}\label{3.30--}
C=\max\left\{\frac{n}{2}\log \frac{(C'+C_1)}{b^2}, \frac{n}{2}\log\frac{C_2}{\tilde{a}^2}\right\}+1,
\end{align}
where $C'=\max_{\overline{B_{2}}}|u_1\Delta^{Ch}_\omega \varphi_1+2\nabla \varphi_1 \cdot\nabla u_1|$. Hence, it is obvious that $u_+=\varphi_1 u_0+C$ is a smooth function on $M$ satisfying $u_+\geq C\geq u_-$ on $M$. 

We claim that $u_+=\varphi_1 u_0+C$ is an upper solution of \eqref{4.8..}. Indeed, if $x\in B_{1}$, then using $\varphi_1|_{B_1}\equiv 1$, \eqref{4.6-}, \eqref{4.11.-} and \eqref{4.14'-}, we obtain
\begin{equation}
\begin{aligned}
-\Delta^{Ch}_\omega u_++k-Ke^{\frac{2}{n}u_+}&=-\Delta^{Ch}_\omega u_0+k-Ke^{\frac{2}{n}(u_0+C)}\\
&=-\Delta^{Ch}_\omega u_1+k-Ke^{\frac{2}{n}( u_1+C)}\\
&\geq 2C_1-C_1-\epsilon e^{\frac{2}{n}( \max_{\overline{B_1}}u_1+C)} \\
&\geq 0,\quad \text{in } B_1,
\end{aligned}
\end{equation}
provided that
\begin{align}
K(x)\leq \epsilon=C_1e^{-\frac{2}{n}(\max_{\overline{B_1}}u_1+C)}.
\end{align}
If $x\in B_2\setminus D$, then using \eqref{4.6-}, \eqref{3.103}, \eqref{4.11.-}, \eqref{3.30--}, $0\leq \varphi_1\leq 1$ and $u_1(x)\geq 0$, we obtain
\begin{equation}
\begin{aligned}
-\Delta^{Ch}_\omega u_++k-Ke^{\frac{2}{n}u_+}&=-u_1\Delta^{Ch}_\omega \varphi_1-2\nabla \varphi_1 \cdot\nabla u_1-\varphi_1\Delta^{Ch}_\omega u_1+k-Ke^{\frac{2}{n}u_+}\\
&=-u_1\Delta^{Ch}_\omega \varphi_1-2\nabla \varphi_1 \cdot\nabla u_1+2C_1\varphi_1+k-Ke^{\frac{2}{n}u_+}\\
&\geq -C'-C_1+b^2e^{\frac{2}{n}\varphi_1 u_1}e^{\frac{2}{n}C}\\
&\geq -C'-C_1+b^2e^{\frac{2}{n}C}\\
&\geq 0,\quad \text{in } B_2\setminus D,
\end{aligned}
\end{equation}
where $C'=\max_{\overline{B_2}}|u_1\Delta^{Ch}_\omega \varphi_1+2\nabla \varphi_1 \cdot\nabla u_1|$. If $x\in M\setminus \text{supp} \varphi_1$, then $u_+(x)=C$. Hence
\begin{equation}
\begin{aligned}
-\Delta^{Ch}_\omega u_++k-Ke^{\frac{2}{n}u_+}\geq -C_1+b^2e^{\frac{2}{n}C}\geq 0,\quad \text{in } M\setminus \text{supp} \varphi_1,
\end{aligned}
\end{equation}
where we have used \eqref{4.6-}, \eqref{3.103} and \eqref{3.30--}. Therefore, $u_+$ is a smooth function on $M$ with $-\Delta^{Ch}_\omega u_++k-Ke^{\frac{2}{n}u_+}\geq 0$ in $M$. Moreover,  by $\varphi_1, u_0\geq 0$ and \eqref{3.30--}, we have $u_+=\varphi_1 u_0+C\geq u_-$, and $\sup_M u_+\leq \|u_1\|_{C^0(\overline{D_2})}+C<+\infty$. By Proposition \ref{prop4.1.}, we deduce that there exists a smooth bounded solution of \eqref{3.104}.
\qed

\section{Nonexistence results of prescribed Chern scalar curvature}
In this section, we will apply the generalized maximum principle (Proposition \ref{prop4.1}) to establish a nonexistence result for prescribed Chern scalar curvature on complete noncompact Hermitian manifolds.
\begin{theorem}[=Theorem \ref{theorem1.9}]
Let $(M^n, \omega)$ be a complete noncompact Hermitian manifold with 
\begin{align}
Ric^{(2)}(X, \overline{X})\geq-C_1(1+r(x))^2,
\end{align}
\begin{align}
S^{Ch}(\omega)\geq 0,
\end{align}
and the torsion 
\begin{align}
\|T^{Ch}(X,Y)\|\leq C_2(1+r(x)),
\end{align}
where $X, Y\in T^{1,0}_{x}M$ satisfy $\|X\|=\|Y\|=1$, $r(x)$ is the Riemannian distance between a fixed point $x_0$ and $x$ in $(M, \omega)$, and $C_1, C_2$ are two positive constants. Let $K$ be a smooth function on $M$ satisfying
\begin{equation}
\begin{aligned}
&K(x)\leq 0, \quad \forall\ x\in M,\\
&K(x)\leq -b^2, \quad\forall\ x\in M\setminus D,
\end{aligned}
\end{equation}
where $D$ is a compact subset of $M$, and $b>0$ is a constant. Then the Hermitian metric $\omega$ cannot be conformally deformed to another Hermitian metric $\tilde{\omega}$ such that its Chern scalar curvature is $K$.
\end{theorem}
\proof We prove it by contradiction. Suppose that there is a function $u\in C^\infty(M)$ such that $\tilde{\omega}=e^{\frac{2}{n}u}\omega$ and the Chern scalar curvature of $(M, \tilde{\omega})$ is $K$. Then $u$ satisfies
\begin{align}\label{5.5}
-\Delta^{Ch}_\omega u+k=Ke^{\frac{2}{n}u},
\end{align}
where $k=S^{Ch}(\omega)\geq 0$. Since $k\geq 0$ and $K\leq (\not\equiv) 0$, we have that $u$ is not a constant and
\begin{align}
\Delta^{Ch}_\omega u\geq k.
\end{align}
Applying the classical maximum principle to the above inequality, we conclude that $u$ cannot attain the maximum on $M$. Set
\begin{align}
v=e^{\frac{2}{n}u}>0.
\end{align}
By \eqref{5.5}, we obtain
\begin{align}\label{5.7}
-\Delta^{Ch}_\omega v+\frac{|\nabla v|^2}{v}+\frac{2}{n}kv=\frac{2}{n}Kv^2.
\end{align}
Since $k\geq 0$, it follows from \eqref{5.7} that
\begin{align}\label{5.9}
\Delta^{Ch}_\omega v\geq -\frac{2}{n}Kv^2.
\end{align}
Let us consider the following smooth function on $M$:
\begin{align}
\phi=\frac{1}{(v+1)^a},
\end{align}
where $0<a<\frac{1}{2}$. Since $\phi>0$, applying Proposition \ref{prop4.1} to $-\phi$, there exists a sequence $\{x_m\}_{m\in\mathbb{N}_+}\subset M$ such that
\begin{align}\label{5.11}
\lim_{m\rightarrow +\infty}\phi(x_m)=\inf_M\phi,
\end{align}
\begin{align}\label{5.12}
\|\nabla \phi(x_m)\|<\frac{1}{m},
\end{align}
\begin{align}\label{5.13}
\Delta^{Ch}_\omega\phi(x_m)>-\frac{1}{m},
\end{align}
for each $m\in \mathbb{N}_+$. A simple computation yields
\begin{align}\label{5.14}
\nabla \phi=-a\frac{\nabla v}{(v+1)^{a+1}},
\end{align}
\begin{align}\label{5.15}
\Delta^{Ch}_\omega\phi=-a\frac{\Delta^{Ch}_\omega v}{(v+1)^{a+1}}+a(a+1)\frac{|\nabla v|^2}{(v+1)^{a+2}}.
\end{align}
By \eqref{5.12} and \eqref{5.14}, we have
\begin{align}\label{5.16}
|\nabla v|^2(x_m)<\frac{1}{a^2m^2}(v(x_m)+1)^{2a+2}.
\end{align}
Substituting \eqref{5.9} and \eqref{5.16} into \eqref{5.15}, it follows from \eqref{5.13} that
\begin{align}
-\frac{1}{m}<\Delta^{Ch}_\omega\phi(x_m)\leq \frac{2a}{n}\cdot\frac{K(x_m)v^2(x_m)}{(v(x_m)+1)^{a+1}}+\frac{a+1}{am^2}(v(x_m)+1)^a,
\end{align}
which implies
\begin{align}
 -\frac{2a}{n}\cdot\frac{K(x_m)v^2(x_m)}{(v(x_m)+1)^{2a+1}}<\frac{1}{m(v(x_m)+1)^a}+\frac{a+1}{am^2}.
\end{align}
According to \eqref{5.11}, we have $\lim_{m\rightarrow +\infty}v(x_m)=\sup_M v$. Since $u$ cannot attain the maximum on $M$, neither can $v$. Hence, $x_m\in M\setminus D$ for sufficiently large $m$. Thus, for $m$ large enough,
\begin{align}\label{5.19}
 -\frac{2a}{n}\cdot\frac{(\sup_{M\setminus D}K)\cdot v^2(x_m)}{(v(x_m)+1)^{2a+1}}<\frac{1}{m(v(x_m)+1)^a}+\frac{a+1}{am^2}.
\end{align}
If $\sup_M v=+\infty$, then taking $m\rightarrow +\infty$ in \eqref{5.19} yields a contradiction, because $\sup_{M\setminus D}K<0$ and $0<a<\frac{1}{2}$. Therefore, $\sup_M v<+\infty$. Now taking $m\rightarrow +\infty$ in \eqref{5.19} again, we get
\begin{align}
\sup_{M\setminus D}K\cdot \sup_Mv^2\geq 0,
\end{align}
which is also a contradiction, since $\sup_{M\setminus D}K<0$. Therefore, there is no function $u\in C^\infty(M)$ solving \eqref{5.5}. Consequently, the Hermitian metric $\omega$ cannot be conformally deformed to another Hermitian metric $\tilde{\omega}$ such that its Chern scalar curvature is $K$.
\qed

\bigskip

Weike Yu

School of Mathematical Sciences, 

Ministry of Education Key Laboratory of NSLSCS,

Nanjing Normal University,

Nanjing, 210023, Jiangsu, P. R. China,

wkyu2018@outlook.com

\bigskip

\end{document}